\documentclass[12pt,reqno]{amsart}

\usepackage[latin1]{inputenc}
\usepackage{amsmath}
\usepackage{amsfonts}
\usepackage{amssymb}
\usepackage{graphics}

\usepackage{enumerate}
\usepackage{amssymb,amsmath,amsthm,amscd}
\usepackage{latexsym,verbatim,graphicx,amsfonts}

\usepackage{amscd}
\usepackage{amsmath}
\usepackage{amssymb}
\usepackage{amsthm}
\usepackage{latexsym}
\usepackage{verbatim}
\usepackage{hyperref}

\theoremstyle{plain}
\newtheorem{theorem}{Theorem}[section]
\newtheorem{thm}[theorem]{Theorem}
\newtheorem{prob}[theorem]{Problem}

\theoremstyle{definition}

\newtheorem{defn}[theorem]{Definition}

\theoremstyle{remark}
\newtheorem{rem}[theorem]{Remark}
\newtheorem*{thm*}{Theorem}

\newcommand{\C}{\mathbb{C}}
\newcommand{\Q}{\mathbb{Q}}

\newcommand{\D}{\mathbb{D}}
\newcommand{\N}{\mathbb{N}}

\newcommand{\McC}{\raise.5ex\hbox{c}}

\title[Dirichlet inner functions]{A characterization of Dirichlet-inner functions}
\author[Seco]{Daniel Seco}
\address{Instituto de Ciencias Matem\'aticas, Calle Nicol\'as Cabrera, 13-15, 28049 Madrid,
Spain.} \email{daniel.seco@icmat.es}

\thanks{
The author thanks Ministerio de Econom\'ia y Competitividad, Spain,
for support through the Severo Ochoa Excellence Grant SEV-2015-0554
and the Research Project MTM2014-52865-P}
\date{\today}

\begin{document}

\maketitle

\begin{abstract}
We study a concept of inner function suited to Dirichlet-type
spaces. We characterize Dirichlet-inner functions as those for which
both the space and multiplier norms are equal to 1.
\end{abstract}

\section{Introduction}\label{Intro}

\subsection{Dirichlet-type and weighted Hardy spaces.} The Hardy space $H^2$ is a well known object in complex
function theory consisting of all holomorphic functions $f$ over the
unit disc $\D$ of the complex plane for which the Taylor
coefficients $\{\hat{f}(k)\}_{k \in \N}$ are in $\ell^2$, with the
corresponding norm defined as
\begin{equation}\|f\|_{H^2}=\left(\sum |\hat{f}(k)|^2\right)^{1/2}.\end{equation} We
refer the reader to \cite{Duren, Garnett} for the classical theory
of $H^2$. $H^2$ admits an integral representation in terms of the
$L^2$ norm with respect to Lebesgue measure over the unit circle.

Consider a measure $\mu$ defined over $\overline{\D}$ and define the
\emph{Dirichlet-type space induced by $\mu$}, $D_{\mu}$, as the
space of holomorphic functions $f$ over $\D$ with finite norm
\begin{equation}\|f\|_{D_{\mu}}:= \left(\|f\|^2_{H^2} + \int_{\D}
|f'(z)|^2 d\mu(z)\right)^{1/2}. \end{equation} We denote the inner
product associated to this norm as $\left< \cdot, \cdot
\right>_{D_\mu}$.

Important particular cases include, again, the Hardy space (for
instance, when $\mu$ is any absolutely continuous measure with
bounded density and compact support), and the classical
\emph{Dirichlet space}, denoted $D$, where $\mu$ is the normalized
Lebesgue measure of area $A$ over $\D$. The theory about the
Dirichlet space has seen an increased interest in the last few
decades. Useful references on Dirichlet-type spaces include
\cite{Arcetal,Primer,Ross}.

Here we focus on the particular case of a radial weight, that is,
the case in which $d\mu(z)=d\mu(|z|)$ for all $z \in \D$, because
such spaces also admit a norm in terms of the Taylor coefficients:
Let $\omega =\{\omega_k\}_{k \in \N}$ be a sequence of positive real
numbers and define the \emph{weighted Hardy space induced by
$\omega$}, $H^2_\omega$, as the space of all holomorphic functions
$f$ over $\D$ with Taylor coefficients $\{\hat{f}(k)\}_{k \in \N}$
with finite norm \begin{equation}\|f\|_{\omega}:=\left(\sum
|\hat{f}(k)|^2 \omega_k \right)^{1/2}.\end{equation} We will
normalize $\omega_0=1$, and assume that there exists some constant
$C>1$ such that for all $k \in \N$ we have
\begin{equation}\label{eqn302} \omega_k \leq C k^2
\end{equation} and
\begin{equation}\label{eqn301} \omega_k \leq \omega_{k+1}.\end{equation} The first property ensures these weighted Hardy spaces
are also (radially weighted) Dirichlet-type spaces while the second
guarantees that the \emph{shift operator} $S$ (of multiplication by
the independent variable $z$) becomes a bounded operator that
increases the norm of every element in the space. We say that a
space satisfying \eqref{eqn301} has the \emph{expansive shift
property}. This excludes the classical Bergman space $A^2$, where
$\omega_k=1/(k+1)$ but it is satisfied by $H^2$ and $D$.
Dirichlet-type spaces with a radial weight and the expansive shift
property will be our setting.

\subsection{Inner functions and multipliers} A function $f \in H^2$ is
called \emph{inner} if it satisfies the following two conditions:
\begin{equation}\label{cond1}|f(z)| \leq 1 \quad \forall z \in \D\end{equation}
\begin{equation}\label{cond2}|f(e^{i \theta})| = 1 \quad a.e.\theta \in [0,
2\pi)\end{equation} Inner functions play a major role in the
function theory of $H^2$, and in particular, in the theory of
invariant subspaces for the shift operator. A space closely related
with $H^2$ is the space $H^\infty$ of bounded analytic functions
over $\D$. A function $g$ is in $H^\infty$ if it has finite norm, as
defined by
\begin{equation}\label{eqn201}\|g\|_{\infty}:= \sup_{z\in \D}
|g(z)|.\end{equation} Since the $L^2$ norm over the circle
represents the norm in $H^2$, it is easy to see that $\|f\|_{H^2}
=1$ for any inner $f$. The starting point for the present article is
the following observation:
\begin{rem} In the definition of inner functions given above we
could express conditions \eqref{cond1} and \eqref{cond2} as
\begin{equation}\label{cond3}\|f\|_{H^2}=\|f\|_{H^\infty}
=1.\end{equation}\end{rem}

Back in the setting of Dirichlet-type spaces $D_{\mu}$, we say that
$h$ is a \emph{multiplier} of $D_{\mu}$ if multiplying by $h$ is a
bounded operator $M_h$ acting on the space $D_{\mu}$. If the norm of
the operator is identified with the norm of $h$ as a multiplier
$\|h\|_{M_{D_\mu}}$, this defines a Banach algebra called the
\emph{space of multipliers} of $D_{\mu}$. As a relevant example,
$H^\infty$ is the space of multipliers of $H^2$ with the norm
\eqref{eqn201}. Observe that the choice of equivalent norm in a
space affects the multiplier norm but it does not affect whether or
not certain function is a multiplier.

When the concept of inner functions has been extended to
Dirichlet-type (and other) spaces, inner functions have been
described in terms of some orthogonality properties of the function.
Denote by $\delta_{\cdot,\cdot}$ the Kronecker delta.
\begin{defn}\label{inner} We say that a function $f$ is $D_\mu$-inner if, for
$j\in \N$,
\begin{equation}\label{cond4}
\left< z^jf, f \right> = \delta_{0,j}.
\end{equation}
\end{defn}
It is easy to check that inner functions are exactly $H^2$-inner
functions. The systematic study of $D_\mu$-inner functions goes back
to the work of Richter in \cite{Richter}, who showed that $D$-inner
functions are also connected with invariant subspaces of the shift
on $D$. However, examples of some $D_\mu$-inner functions can
already be found in the work of Shapiro and Shields (\cite{ShaShi}).
They provide examples of inner functions with any prescribed finite
zero set on the unit disc. For further background on the relation
between generalized inner functions, invariant subspaces, and other
topics, see \cite{BFKSS}.

In 1992, Richter and Sundberg (\cite{RichSund92}) proved the
following connection between inner functions and multipliers, as
part of a stronger result:
\begin{thm}\label{RS92}
Let $\mu$ be a harmonic weight and $f$, $D_\mu$-inner. Then
\begin{equation}\label{eqn202} \|f\|_{M_{D_\mu}} =1.
\end{equation}
\end{thm}
Another proof was given by Aleman in \cite{Aleman}. See also
\cite{Primer}, Theorem 8.3.9.

\section{Main result}

Our main goal in this paper is to prove a reciprocal to Theorem
\ref{RS92}. As a byproduct, we obtain a relatively small set of test
functions that contain all the information on whether a given
function is inner.

\begin{thm}\label{main}
Let $f \in D_\mu=H^2_\omega$. Then $f$ is $D_\mu$-inner whenever
\begin{itemize}
\item[(a)] $\|f\|_{D_\mu}=1$ and for all $k \in \N \backslash \{0\}$ and all $\lambda \in \C$,
$$\|fg_{k, \lambda}\|^2_{D_\mu} \leq \omega_k + |\lambda|^2,$$ where
$g_{k, \lambda}(z) = z^k+\lambda$,\end{itemize} and this holds true
whenever
\begin{itemize}
\item[(b)] $$\|f\|_{D_\mu}=\|f\|_{M_{D_\mu}}=1.$$
\end{itemize}
\end{thm}

Given Theorem \ref{RS92}, this is a characterization of
Dirichlet-inner functions.

\begin{proof} Clearly, (b) implies (a),
since (a) only requires that
$$\|fg\|_{D_\mu} \leq \|g\|_{D_\mu}$$ holds for a subset of
all possible $g \in D_\mu$. Therefore what remains is to show that
(a) implies $D_\mu$-inner.

To see this, assume $f$ is not $D_\mu-$inner but $\|f\|_{D_\mu}=1$
and let $k \in \N \backslash \{0\}$ such that
\begin{equation}\label{eqn501}
\left< z^k f, f \right> \neq 0.
\end{equation}
If we find a value of $\lambda$ such that
\begin{equation}\label{eqn502}
\|(z^k+\lambda)f\|^2_{D_\mu} > \|(z^k + \lambda)\|^2_{D_\mu}=
\omega_k + |\lambda|^2,
\end{equation} then (a) is not satisfied and we will be done.

Decompose the right-hand side of \eqref{eqn502} in terms of inner
products as
\begin{equation}\label{eqn503}
\|(z^k+\lambda)f\|^2_{D_\mu} = \|z^k f\|^2_{D_\mu} + |\lambda|^2
\|f\|^2_{D_\mu} + 2 Re (\overline{\lambda} \left<z^k f, f\right>).
\end{equation}

Since the shift is expansive, we have that
\begin{equation}\label{eqn504}\|z^k f\|^2_{D_\mu} \geq
\|f\|^2_{D_\mu}=1.
\end{equation}
Using that $\|f\|_{D_\mu}=1$ and \eqref{eqn501}, we may choose
\begin{equation}\label{eqn505}\lambda =\frac{\left<z^k f, f\right>}{|\left<z^k f, f\right>|^2}
\frac{\omega_k}{2},
\end{equation} and this will yield
\begin{equation}\label{eqn506}
\|(z^k+\lambda)f\|^2_{D_\mu} \geq 1 + |\lambda|^2 + \omega_k.
\end{equation} That is, this choice of $\lambda$ achieves
\eqref{eqn502} completing the proof.
\end{proof}

\section{Further remarks}

\begin{itemize}
\item[(A)] The multipliers of $D$ were characterized by Stegenga in \cite{Stegenga} in terms of Carleson measures. See also
\cite{Arcetal}. It would be interesting to describe properties (a)
or (b) in Theorem \ref{main} from Carleson measures.
\item[(B)] In spaces without the expansive shift property, Theorem
\ref{main} fails to hold in general as $D_\mu$-inner does not always
imply (b). For instance, take the classical Bergman space
$A^2=H^2_\omega$ for $\omega_k=1/(k+1)$, which admits a
representation as the $L^2$ integral over the disc. Its space of
multipliers is $H^\infty$ (as for $H^2$). Any nonconstant function
bounded by 1 will have Bergman norm strictly less than 1, and hence
no nonconstant functions satisfy (b), while there is a plethora of
nonconstant $A^2$-inner functions (as made clear in any of
\cite{ARS, Heden, DKSS, BFKSS}). However, our proof of (b)
$\Rightarrow$ (a) $\Rightarrow$ $D_\mu$-inner works in any weighted
Hardy space with the expansive shift property: we do not make use of
the assumption \eqref{eqn302}. Therefore, it seems natural to ask
whether Theorem \ref{RS92} holds true on all weighted Hardy spaces
with the expansive shift property. It even makes sense to ask for an
analogue to Theorem \ref{main} on reproducing kernel Hilbert spaces
where the shift increases the norm. Advances on the understanding of
the multiplier properties of general spaces inner functions have
recently been made in \cite{Aleetal} and other references therein.
This may yield an additional way to prove the result in this article
in a more general environment. The dichotomic situation between
spaces with or without the expansive shift property is also present
in previous work (\cite{Zeros}). In spaces where the shift is
contractive, an analogue proof to that of our theorem will show that
contractive divisors are inner, although in the most classical case
this is already well known.
\item[(C)] Condition (a) in Theorem \ref{main} tells us that we only
need to check the multiplication properties on a sequence of
unidimensional vector subspaces of $D_\mu$. In fact, it can be shown
that it suffices to test on a countable set of functions (for
example, by taking $\lambda$ on $\Q + i \Q$). If moreover, we are
performing the test on a function whose Taylor coefficients are
real, then $\lambda$ can be taken to be real.
\item[(D)] From Theorem \ref{main} we learn that inner functions are exactly those
functions that solve the extremal problem
\begin{equation}\label{eqn401} \inf \left\{
\frac{\|f\|_{M_{D_\mu}}}{\|f\|_{D_\mu}}: f \in D_\mu
\backslash\{0\}\right\}.
\end{equation}
It seems conceivable that a similar property is satisfied in other
spaces like the Bergman space $A^2$. We conclude by proposing a few
problems along this line:
\begin{prob} Let $M$ be a $z-$invariant subspace of $A^2$ containing at least one bounded function. Let $g$
be a function solving the extremal problem
\begin{equation}\label{eqn402} \inf \left\{
\frac{\|f\|_{H^\infty}}{\|f\|_{A^2}}: f \in M\backslash\{0\}
\right\}.
\end{equation} Is it always true that $g$ is a constant multiple of an $A^2-$inner function?
\end{prob}
 In \cite{BH}, the authors showed that
there exist invariant subspaces of $A^2$ containing no bounded
functions and, therefore, we need the hypothesis on bounded
functions. See also \cite{AR1} and the references therein.

A follow up question is a reciprocal to this:
\begin{prob}
Does the orthogonal projection of the constant function 1 onto $M$
solve the extremal problem \eqref{eqn402}?
\end{prob}
The same questions make sense in any reproducing kernel Hilbert
space where inner functions are defined.
\end{itemize}

\noindent\textbf{Acknowledgements.}  The author would like to thank
A. Borichev, M. Hartz and T. Le for their useful comments.

\end{document}